\documentclass[12pt]{article}

\def\@eqnnum{\rm ([section].\theequation)}
\usepackage[left=3cm,top=4cm,right=3cm,bottom=2cm,nohead,foot=1.25cm]{geometry}
\usepackage{amsthm}
\usepackage{amsfonts}
\usepackage{pictexwd,dcpic}
\usepackage{amsmath}
\usepackage{amssymb}
\numberwithin{equation}{section}

\DeclareMathOperator{\fun}{\hspace{1mm}is\hspace{1mm}a\hspace{1mm}
constant\hspace{1mm} function}
\DeclareMathOperator{\for}{\hspace{1mm}for\hspace{1mm}}
\DeclareMathOperator{\ifff} {iff}
\DeclareMathOperator{\iif}{if}
\newtheorem{df}{Definition}[section]
\newtheorem{p}[df]{Proposition}
\newtheorem{corol}[df]{Corollary}
\date{}
\title{\bf  An algebraic analysis framework \\for quantum calculus}
\author{Piotr Multarzy\'nski \\{\it Faculty of Mathematics and Information Science}\\
 {Warsaw University of Technology} \\
{\it 00-661 Warsaw, Pl. Politechniki 1, Poland}\\
e-mail: multarz@mini.pw.edu.pl}
\begin{document}
\maketitle
\begin{abstract}
An algebraic analysis framework for quantum calculus is proposed.
The quantum derivative operator $D_{\tau ,\sigma}$ is based on two
commuting bijections $\tau$ and $\sigma$ defined on an arbitrary set
$M$ equipped with a tension structure determined by a single tension
function $\theta$, i.e. a 1-dimensional case is analyzed here. The
well known cases, i.e. $h$- and $q$-calculi together with their
symmetric versions, can be obtained owing to special choice of
mappings $\tau$ and $\sigma$.
\end{abstract}
{\bf Keywords:} Quantum calculus;  difference quotient operator; right invertible operator \\
{\bf MSC 2000: 12H10, 39A12, 39A70}
\section{Introduction}
The term "algebraic analysis" is used by many authors to indicate an
algebraic approach to analytic problems and, in fact, it is used in
many different senses. In the present paper this term we use in the
sense of D. Przeworska-Rolewicz \cite{DPR} since our main interest
here is the calculus of right invertible operators. The examples of
such operators can be the usual derivative $\frac{d}{dx}$ as well as
the divided difference operators studied in quantum calculus
\cite{kac}. In Section \ref{qintegration} we interpret the quantum
$h$- and $q$-definite integrals within the algebraic analysis
framework. Then, in Section \ref{tausigmacalculus} some more general
setting is proposed. Namely, we analyze linear operators defined on
function (commutative) algebras and satisfying certain product rules
being the modified versions of the Leibniz rule. Such operators have
many properties which are quite analogous with the corresponding
ones for differential operators. In parallel, there is a natural
possibility to define some kind of algebraic integration associated
with right invertible operators. The algebraic concept of definite
integration, with respect to a given right invertible linear
operator, has been defined by using initial operators within the
algebraic analysis framework \cite{DPR}.

For the reader's convenience, below we present some basics of
algebraic analysis and quantum calculus.
\section{Calculus of right invertible operators}
Quantum calculus is in fact a sort of a discrete  calculus, in which
some discrete versions of differentiation and integration are
studied. In the present paper we are going to compare the proposal
of quantum calculus integration  with the corresponding general idea
of integration based on the calculus of right invertible operators.
For the reader's convenience we give a short survey of the basic
concepts concerning the right invertible operators but the
comprehensive treatment of the topic one can find in Reference
\cite{DPR}.

Let $X$ be a linear space over $\mathbb R$ and $L(X)$ be the family
of all linear operators in $X$ with the domains being linear
subspaces of $X$. Then, for any $A\in L(X)$, let ${\cal D}_{A}$
denote the domain of $A$ and let $L_{0}(X)=\{A\in L(X): {\cal
D}_{A}=X\}$. By the space of constants of an operator $D\in L(X)$ we
shall mean the set $Z_{D}=ker D$. A linear operator $D\in L(X)$ is
said to be right invertible if $DR=I$, for some linear operator
$R\in L_{0}(X)$ called a right inverse of $D$ and $I=id_{X}$. The
family of all right invertible operators in $X$ will be denoted by
${\cal R}(X)$. In turn, by ${\cal
R}_{D}=\{R_{\gamma}\}_{\gamma\in\Gamma}$ we shall denote the family
of all right inverses of a given $D\in {\cal R}(X)$. If $R\in{\cal
R}_{D}$ is a given right inverse of $D\in{\cal R}(X)$, the family
${\cal R}_{D}$ is characterized by
\begin{equation}\label{RA}
{\cal R}_{D}=\{R+(I-RD)A: A\in L_{0}(X)\}.
\end{equation}
Consider a family of right invertible operators $D_i\in{\cal R}(X)$
and a corresponding family of their right inverses $R_i\in{\cal
R}_{D}$, for $i=1,\ldots ,n$ and some $n\in\mathbb N$. Then, the
composition $D=D_1\ldots D_n$ is right invertible, i.e. $D\in{\cal
R}(X)$, and one of its right inverses $R\in{\cal R}_{D}$ is given by
\begin{equation}\label{composeinverse}
R=R_n\ldots R_1\; .
\end{equation}
For any $x, y\in X$, we say that $y$ is a primitive element of $x$
whenever $Dy=x$. Thus, the element $Rx$ is a primitive element of
$x$, for any $x\in X$ and $R\in {\cal R}_{D}$. The set
\begin{equation}\label{indefint}
{\cal I}(x)=\{y\in X : Dy=x\}
\end{equation}
is called the indefinite integral of a given $x\in X$.
\hspace{1mm}One can easily check, that
\begin{equation}
{\cal R}_{D}x=\{ Rx+(I-RD)Ax: A\in L_{0}(X)\}= {\cal R}_{D}x +
Z_{D}=Rx+Z_{D},
\end{equation}
for any  $R\in {\cal R}_{D}$ and any non-zero element $x\in X$.
Hence, we obtain
\begin{equation}\label{indefintRxZ}
{\cal I}(x)={\cal R}_{D}x + Z_{D}=Rx+Z_{D},
\end{equation}
for any $x\in X$ and $R\in {\cal R}_{D}$.

Any projection operator $F\in L_{0}(X)$ onto $Z_D$, i.e. $F^2=F$ and
$Im F=Z_D$, is said to be an initial operator induced by $D\in{\cal
R}(X)$ and the family of all such operators we denote by ${\cal
F}_D$. For an initial operator $F$ and $x\in X$, the element $Fx\in
Z_D$ is called the initial value of $x$. Additionally, we say that
an initial operator $F\in{\cal F}_D$ corresponds to $R\in {\cal
R}_{D}$ if $FR=0$ or equivalently if
\begin{equation}\label{FR}
F=I-RD,
\end{equation}
on the domain of $D$.  The two families ${\cal R}_D$ and ${\cal
F}_D$ uniquely determine each other. Indeed, for any $R\in{\cal
R}_D$ we define $F=I-RD\in{\cal F}_D$. On the other hand, for any
$F\in{\cal F}_D$, we define $R=R_1-FR_1$, where $R_1\in{\cal R}_D$
can be any since the result is independent of the choice of $R_1\,$.
Thus, for any $\gamma\in\Gamma$ we have $F_\gamma =I-R_\gamma D$ and
consequently ${\cal F}_D=\{F_\gamma\}_{\gamma\in\Gamma}\,$.

By a simple calculation one can verify that $F_\alpha F_\beta
=F_\beta$ and $F_\beta R_\alpha = R_\alpha - R_\beta$, for any
$\alpha , \beta\in\Gamma$. Hence, for any indices $\alpha, \beta,
\gamma\in\Gamma$, there is
\begin{equation}\label{FRgamma}
F_\beta R_\gamma - F_\alpha R_\gamma = F_\beta R_\alpha,
\end{equation}
which means that in fact the left side of equation (\ref{FRgamma})
is independent of $\gamma$. The last property allows one to define
the following definite integration operator
\begin{equation}\label{definiteI}
{\cal I}^{\beta}_{\alpha}=F_\beta R_\gamma - F_\alpha R_\gamma ,
\end{equation}
for any $\alpha, \beta, \gamma\in\Gamma$. Amongst many properties of
the operator ${\cal I}^{\beta}_{\alpha}$ we can mention the most
intuitive one, namely
\begin{equation}
{\cal I}^{\beta}_{\alpha}D=F_\beta - F_\alpha\; .
\end{equation}
 Hence, for any $x\in X$ and its arbitrary
primitive element $y\in X$, i.e. $Dy =x$, we obtain
\begin{equation}
{\cal I}^{\beta}_{\alpha}x=F_\beta y - F_\alpha y\; ,
\end{equation}
 which is called the definite integral of $x$.

To intuitively demonstrate the basic concepts of algebraic analysis,
we end this section with two important examples. In the first
example we take the usual derivative operator $D=\frac{d}{dx}$ while
in the second example we consider $D=D_h$ being the following
difference operator defined by
\begin{equation}\label{d_h}
D_hf(x)=\frac{f(x+h)-f(x)}{h}\; ,
\end{equation}
and giving rise to h-quantum calculus.\vspace{2mm}\\
{\bf Example 1.1} Assume the linear space $X=C^0({\mathbb R})$ (all 
continuous real functions) and $D=\frac{d}{dx}$. Then we recognize
the domain ${\cal D}_D=C^1({\mathbb R})$ (all real functions having
continuous derivative) and the set (linear subspace) of all
constants of $D$ is $Z_D=\{f\in X: f\fun\}$. Since $Z_D$ is a
$1$-dimensional linear space over $\mathbb R$, we shall assume the
identification $Z_D\equiv\mathbb R$. Thus, the initial operators $F$
in this example are projections of $X$ onto $\mathbb R$. To see why
the name "initial operator" is intuitively consonant, it is enough
to notice that $F_a: X\ni f\mapsto f(a)\in Z_D\equiv\mathbb R$,
$a\in\mathbb R$, is a projection operator associating with any $f$
its value at $a$. Obviously, $\{F_a: a\in{\mathbb R}\}\subset{\cal
F}_D$ but $\{F_a: a\in{\mathbb R}\}\neq{\cal F}_D$. The reason is
that any convex combination of initial operators is again an initial
operator \cite{DPR}. For example, one can easily check that
$F_{ab}=\frac{1}{2}(F_a+F_b)\in{\cal F}_D$ and $F_{ab}\neq F_c$, if
$a\neq b$, for any $a,b,c\in\mathbb R$. Therefore, although ${\cal
F}_D$ can be viewed as an indexed family, i.e. ${\cal
F}_D=\{F_\gamma\}_{\gamma\in\Gamma}$, we cannot naturally identify
$\Gamma$ with $\mathbb R$. As an example of a right inverse of
$D=\frac{d}{dx}$ we can take $R: X\rightarrow C^1({\mathbb R})$,
such that $R_af(x)=\int\limits_{a}^{x}f(t)dt$, for a fixed
$a\in\mathbb R$. Let us notice that $F_a$ is the initial
operator corresponding to $R_a$, for any $a\in\mathbb R$.\vspace{2mm}\\
{\bf Example 1.2} Let 
$X={\mathbb R}^{\mathbb R}$ be the linear space of all real
functions and consider $D$ defined by formula (\ref{d_h}), for a
fixed $h > 0$. Evidently, the linear space $Z_D$ consists of all
$h$-periodic functions. Then, the operator $R$ defined by
\begin{equation}\label{R}
Rf(x)=\left\{
\begin{array}{ccl}
-h\sum\limits_{m=0}^{-\left\lfloor\frac{x}{h}\right\rfloor-1}f(x+mh)& & x\in(-\infty , 0)\vspace{1mm}\\
0 &\for &x\in [0, h)\\
h\sum\limits_{m=1}^{\left\lfloor\frac{x}{h}\right\rfloor}f(x-mh)&& x\in [h , \infty)\\
\end{array}
\right.
\end{equation}
fulfils the condition $DR=I$ and therefore it is a right inverse of
$D$. In the above formula the floor brackets $\lfloor\cdot\rfloor$
stand for the integer value function of its argument. Then, let us
define the operator $F$ by formula
\begin{equation}\label{FRexplicite}
Ff(x)=f(x-\left\lfloor\frac{x}{h}\right\rfloor h),
\end{equation}
for any $x\in\mathbb R$. Since $Ff(x+h)=Ff(x)$, for any $x\in\mathbb
R$, the function $Ff$ is $h$-periodic, i.e. $Ff\in Z_D$. On the
other hand, for any function $f\in Z_D$, there is $Ff=f$. Hence, the
operator $F$ is a projection of $X$ onto $Z_D$ and therefore it is
an initial operator induced by $D=D_h$. Moreover, one can check the
property (\ref{FR}) which means that $F$ is the initial operator
corresponding to $R$ given by (\ref{R}). By using formula
(\ref{RA}), the family ${\cal R}_D$ is fully determined by the above
operator $R$, then with the help of (\ref{FR}) we obtain the family
${\cal F}_D$.
\section{Quantum h- and q-calculus}\label{hqcalculus}
In this section we briefly recall the main elements of quantum
calculus but more detailed study of the topic, motivation and many
properties reflecting the analogies with the usual differential
calculus the reader will find in \cite{kac}. For the history of
q-calculus, its relation to other mathematical and physical areas
and the imposing list of references we recommend \cite{Ernst}.

Let $f:{\mathbb R}\rightarrow{\mathbb R}$ be an arbitrary function
and consider the well known difference quotient
\begin{equation}\label{diffquot}
\frac{f(x)-f(x_0)}{x-x_0}\; ,
\end{equation}
for some $x\neq x_0$.
 The limit of the last
expression when $x\rightarrow x_0$, if it exists, defines the
derivative of $f$ at $x_0$. Now, if we take $x=x_0+h$ for a fixed
$h\neq 0$ or $x=qx_0$ for a fixed $q\neq 1$ and do not take the
corresponding limit, we enter the so-called quantum $h$- or
$q$-calculus. For any $f:{\mathbb R}\rightarrow{\mathbb R}$ one
defines its $h$-differential $d_hf$
\begin{equation}\label{h-diff}
d_hf(x)=f(x+h)-f(x),
\end{equation}
and its $q$-differential $\delta_qf$
\begin{equation}\label{q-diff}
\delta_qf(x)=f(qx)-f(x).
\end{equation}
In particular, for the identity mapping $e$ defined on ${\mathbb
R}$, i.e. $e(x)\equiv x$, we have $d_h e(x)=h$ and $\delta_q
e(x)=(q-1)x$. Quite commonly the simplified notation is used, i.e.
$d_h e(x)\equiv d_hx$ and $\delta_q e(x)=\delta_qx$. In
applications, the two versions of quantum calculus (i.e. $h$- or
$q$-calculus) are considered separately, which allows one to denote
both differentials by the same symbol, i.e.  one can write $d_h$ or
$d_q$ (instead $\delta_q$) and recognize them from context. The
above two symbols $d_h$, $\delta_q$ can be viewed as the linear
operators $d_h: f\mapsto d_h f$ and $\delta_q: f\mapsto \delta_q f$
defined on some $\mathbb R$-algebra $\cal A$ of real functions.
However, the algebra $\cal A$ should be invariant w.r.t. the h- or
q-shifts, i.e. functions $x\mapsto f(x+h)$ or $x\mapsto f(qx)$
should be the elements of $\cal A$, for any $f\in\cal A$.

One can easily verify the following Leibniz product rules
\begin{equation}\label{h-Leibniz}
d_h(fg)(x)=d_h(f)(x)g(x+h)+f(x)d_h(g)(x),
\end{equation}
and similarly
\begin{equation}\label{q-Leibniz}
\delta_q(fg)(x)=\delta_q(f)(x)g(qx)+f(x)\delta_q(g)(x).
\end{equation}
The above Leibniz formulae define the corresponding classes ($\cal
A$-modules) of linear difference-like operators, defined on some
$\mathbb R$-algebra ${\cal A}$ of functions.

 Evidently, the above product rules are also fulfilled by operators
$D_h$ and $\Delta_q$, called the quantum derivatives and defined as
\begin{equation}\label{h-deriv}
D_h(f)(x)=\frac{d_hf(x)}{d_he (x)}\equiv \frac{d_hf(x)}{d_hx}
\end{equation}
and similarly
\begin{equation}\label{q-deriv}
\Delta_q(f)(x)=\frac{\delta_qf(x)}{\delta_qe(x)}\equiv\frac{\delta_qf(x)}{\delta_qx}\;
.
\end{equation}
{\bf Remark:} Since $\delta_qe (0)=0$, the expression
$\Delta_q(f)(x)$ is not defined at $x=0$ unless $f'(0)$ does exist.
Therefore, the $q$-calculus can be developed in algebras $\cal A$ of
functions defined on ${\mathbb R}\setminus \{0\}$ or in algebras of
functions defined on ${\mathbb R}$ and differentiable at $x=0$.

In $h$-calculus, an $h$-antiderivative of a function $f: {\mathbb
R}\rightarrow{\mathbb R}$ is defined to be any function $g: {\mathbb
R}\rightarrow{\mathbb R}$ such that $D_hg(x)=f(x)$, for any $x\in
\mathbb R$. The family of all $h$-antiderivatives of a given
function $f$ is called the indefinite $h$-integral and is denoted by
\begin{equation}\label{h-indefintegral}
\int f(x)d_hx\; .
\end{equation}
Then, the definite $h$-integral is defined by formula
\begin{equation}\label{h-defintegral}
\int\limits_{a}^{b} f(x)d_hx= \left\{
\begin{array}{ccl}
h\left(f(a)+f(a+h)+\ldots +f(b-h)\right)&if & a<b \vspace{1mm}\\
0 &if&a=b\; ,\\
-h\left(f(b)+f(b+h)+\ldots +f(a-h)\right)&if& a>b
\end{array}
\right.
\end{equation}
for any $a, b\in\mathbb R$, such that $a$ and $b$ differ by an
integer multiple of $h$.

Concerning $q$-calculus, in this paper we shall assume $q\in (0,
1)\cup (1, \infty )$. This restriction follows from the physical
motivation that the two quantum parameters are usually related by
$q=e^h$. The last exponential relation transforms the real line
$\mathbb R$ onto ${\mathbb R}_+=(0, \infty )$. Consequently, the
$h$-calculus for functions defined on $\mathbb R$ and the
$q$-calculus for functions defined on ${\mathbb R}_+$ can be unified
within a more general framework (generalized quantum calculus). A
$q$-antiderivative of a function $f$ is said to be any function $g$
such that $\Delta_qg(x)=f(x)$. A special $q$-antiderivative, the
so-called Jackson integral of $f$, is formally derived in \cite{kac}
as the geometric series expansion
\begin{equation}\label{Jackson}
g(x)=(1-q)x\sum\limits_{m=0}^{\infty}q^mf(xq^m).
\end{equation}
Then, formula (\ref{Jackson}) is used to define
\begin{equation}
\int\limits_{0}^{b}f(x)d_qx=(1-q)b\sum\limits_{m=0}^{\infty}q^mf(bq^m),
\end{equation}
and finally define  the definite $q$-integral \cite{Jackson}
\begin{equation}\label{qdefinite}
\int\limits_{a}^{b}f(x)d_qx=\int\limits_{0}^{b}f(x)d_qx-\int\limits_{0}^{a}f(x)d_qx\,
,
\end{equation}
 for $0<a<b$.
 Since formula (\ref{Jackson}) has been derived formally, one
needs to examine the conditions when it converges to a
$q$-antiderivative. Within the algebraic analysis framework, used in
this paper, we construct many $q$-antiderivatives which are finite
sums and no condition has to be examined to justify their
convergency. However, the above Jackson integral can be recovered in
this approach provided the corresponding infinite expansion is
convergent.

At the end of this section let us briefly discuss the lack of
symmetry one can notice concerning the product rules
(\ref{h-Leibniz}), (\ref{q-Leibniz}).

On the strength of formulae (\ref{h-Leibniz}), (\ref{q-Leibniz}), 
for any $a,b\in\mathbb R$ such that $a+b=1$, one can write the
following combinations
\begin{equation}\label{h-symalgrule}
d_h(fg)(x)=d_h(f)(x)(ag(x)+bg(x+h))+(bf(x)+af(x+h))d_h(g)(x),
\end{equation}
and analogously
\begin{equation}\label{q-symalgrule}
\delta_q(fg)(x)=\delta_q(f)(x)(ag(x)+bg(xq))+(bf(x)+af(xq))\delta_q(g)(x).
\end{equation}
If $a\neq b$, the above combined formulae (\ref{h-symalgrule}),
(\ref{q-symalgrule}) are equivalent with (\ref{h-Leibniz}) and
(\ref{q-Leibniz}), respectively, i.e. the corresponding classes of
operators defined coincide. On the other hand, for $a=b=\frac{1}{2}$
the corresponding symmetric product rule defines a larger class of
operators, in general. However, for some algebras the symmetric rule
can be equivalent with its all non-symmetric counterparts. A
non-trivial example of an algebra, for which the symmetric product
rule implies all the other ones, is the algebra of polynomials
${\cal A}={\mathbb R}[x]$. Indeed, assume $h=1$, $a=b=\frac{1}{2}$
and consider the symmetric product rule
\begin{equation}\label{h=1-symalgrule}
D(fg)(x)=D(f)(x)\cdot\frac{g(x)+g(x+1)}{2}+\frac{f(x)+f(x+1)}{2}\cdot
D(g)(x)\; ,
\end{equation}
for any $f,g\in{\mathbb R}[x]$.  One can prove that $D=D(1_{\cal
A})\cdot d_1$, where $d_1f(x)=f(x+1)-f(x)$. Hence, any operator $D$
is proportional to the usual difference operator $d_1$ and
consequently it fulfills (non-symmetric) formula (\ref{h-Leibniz}).
In turn, an algebra for which the symmetric product rule
(\ref{h=1-symalgrule}) is weaker than any non-symmetric one is for
example the $\mathbb R$-algebra (of real functions) ${\cal A}=gen
(\{e,z\})$  generated by the identity $e(x)=x$ and the integer
valued function $z(x)=\lfloor x\rfloor$.

\section{An algebraic analysis  approach \newline to quantum integration}\label{qintegration}
In this section we present an  approach to quantum integration
within the algebraic analysis framework \cite{DPR}.

Let us construct the following right inverses $R_{hs}$, $s\in\mathbb
R$. Namely, for  $h<0$
\begin{equation}\label{Rhs-}
R_{hs}f(x)=\left\{
\begin{array}{ccl}
\sum\limits_{m=1}^{-\lfloor
\frac{x-s}{-h}\rfloor}hf(x-mh)& & x\in(-\infty , s)\vspace{1mm}\\
0 &\for &x\in [s, s-h)\\
-\sum\limits_{m=0}^{\lfloor\frac{x-s}{-h}\rfloor -1}hf(x+mh)&& x\in [s-h , \infty)\\
\end{array}
\right.
\end{equation}
 and for  $h>0$
\begin{equation}\label{Rhs+}
R_{hs}f(x)=\left\{
\begin{array}{ccl}
-\sum\limits_{m=0}^{-\lfloor\frac{x-s}{h}\rfloor -1}hf(x+mh)& & x\in(-\infty , s)\vspace{1mm}\\
0 &\for &x\in [s, s+h)\\
\sum\limits_{m=1}^{\lfloor\frac{x-s}{h}\rfloor }hf(x-mh)& & x\in [s+h , \infty)\\
\end{array}
\right.
\end{equation}
By a straightforward calculation one can show that $D_hR_{hs}=I$,
for any $h\neq 0$ and $s\in\mathbb R$.

Then, on the strength of formula (\ref{indefintRxZ}),
the indefinite $h$-integral of a function $f\in{\mathbb R}^{\mathbb R}$
can be written as
\begin{equation}
\int f(x)d_hx=R_{hs}f(x)+Z_{D_h}\equiv R_{h0}f(x)+Z_{D_h},
\end{equation}
where $s\in\mathbb R$ is an arbitrarily fixed index, e.g. $s=0$, and
the notation (\ref{h-indefintegral}) was used.

Define the operators $F_{hs}$ by
\begin{equation}\label{Fa}
F_{hs}f(x)=f(x-\left\lfloor\frac{x-s}{|h|}\right\rfloor\cdot |h|),
\end{equation}
for any $s\in\mathbb R$.

 One can verify that  $F_{hs}$ are the initial operators induced by $D_h$, for any $s\in\mathbb R$.
 Indeed, the function $F_{hs}f$ is $h$-periodic, since
$$F_{hs}f(x+h)=f(x+h-\left\lfloor\frac{x+h-s}{|h|}\right\rfloor\cdot |h|)=$$
$$=f(x+h-\left(\left\lfloor\frac{x-s}{|h|}\right\rfloor+\frac{h}{|h|}\right)\cdot |h|)=F_{hs}f(x).$$ Moreover,
for any $h$-periodic function $f$, we have the evident identity
$F_{hs}f=f$, which proves that $F_{hs}$ is a surjective projection
onto $Z_{D_{h}}$. Therefore $F_{hs}$ is an initial operator induced
by $D_h$, for any $s\in\mathbb R$. One can also verify that the
initial operators $F_{hs}$ correspond to (\ref{Rhs-}) and
(\ref{Rhs+}), respectively.

In turn, the definite $h$-integrals are defined in a general manner
by using formula (\ref{definiteI}). Within this approach we obtain a
large class of definite $h$-integrals, with the integration limits
being arbitrary (indices of) initial operators. Below we consider
definite $h$-integrals determined by the (particular) initial
operators $F_{hs}$, for $s\in\mathbb R$. As we shall see, these
integrals can be used to obtain the ordinary $h$-definite integrals
defined by formula (\ref{h-defintegral}). Namely, by formula
(\ref{definiteI}), for any $a,b\in\mathbb R$ and the corresponding
(particular) initial operators $F_{ha}$,  $F_{hb}$ we obtain

\begin{equation}\label{Iab}
{\cal I}^{b}_{a}=F_{hb}R_{h0}-F_{ha}R_{h0}.
\end{equation}
The concrete right inverse $R_{h0}$, for $s=0$, is used above only
for simplicity since the result is independent of this choice.
Assume $h>0$ and calculate the definite integral of a function $f$
at $x$, i.e.
\begin{equation}
{\cal
I}^{b}_{a}f(x)=R_{h0}f(x-\left\lfloor\frac{x-b}{h}\right\rfloor\cdot
h)-R_{h0}f(x-\left\lfloor\frac{x-a}{h}\right\rfloor\cdot h).
\end{equation}
Then, the ordinary $h$-definite integrals, defined intuitively in
\cite{kac}, are obtained here as the value of ${\cal I}^{b}_{a}f(x)$
at any point $x=a+kh$, $k\in\mathbb Z$. Indeed, assume
$0<k\in\mathbb Z$ and $b=a+kh$. Then, for $x=a$ we obtain
$$
{\cal I}^{b}_{a}f(a)=R_{h0}f(a+kh)-R_{h0}f(a)=
$$
$$
=\sum_{j=0}^{k-1}(R_{h0}f(a+(j+1)h)-R_{h0}f(a+jh))=
\sum_{j=0}^{k-1}hD_hR_{h0}f(a+jh)=
$$
$$
=\sum_{j=0}^{k-1}hf(a+jh)=h(f(a)+f(a+h)+\ldots +f(b-h)).
$$

If $a=b$, the result is obviously ${\cal I}^{b}_{a}f(a)=0$. In turn, for $a>b$
and $a=b+kh$, for some $0<k\in\mathbb Z$,
we have
$$
{\cal I}^{b}_{a}f(b)=R_{h0}f(b)-R_{h0}f(b+kh)=
$$
$$
=\sum_{j=0}^{k-1}(R_{h0}f(b+jh)-R_{h0}f(b+(j+1)h))=
-\sum_{j=0}^{k-1}hD_hR_{h0}f(b+jh)=
$$
$$
-\sum_{j=0}^{k-1}hf(b+jh)=-h(f(b)+f(b+h)+\ldots +f(a-h)).
$$
The above calculation demonstrates how the ordinary $h$-definite
integrals, defined by (\ref{h-defintegral}), emerge from the
algebraic analysis approach used here.

  Directly from definition
of the initial operator concept and from (\ref{Iab}), we conclude
that ${\cal I}^{b}_{a}f\in Z_D$, i.e. it is an $h$-periodic function
and ${\cal I}^{b}_{a}f(x)={\cal I}^{b}_{a}f(a)$, for any $x\in
a+h\mathbb Z$. Let us emphasize the conceptual difference between
definitions (\ref{h-defintegral}) and (\ref{Iab}). In $h$-calculus,
by formula (\ref{h-defintegral}) one defines the definite integral
to be a scalar-valued linear functional, while in the algebraic
analysis approach the corresponding definite integral value is an
$h$-periodic function (non-constant, in general). The above two
formulations of definite integrals are equivalent for functions
defined on the domain $a+h\mathbb Z$, for some fixed $a, h\in\mathbb
R$. \\
{\bf Remark:} Imagine that an action functional of a physical
system is defined as an $h$-integral of some lagrangian.
Consequently, such an action is $h$-periodic and its $h$-periodicity
can be viewed as a physical symmetry giving rise to a corresponding
conservation law.

Concerning $q$-calculus, we shall work here with functions $f$
defined on the domain $(0, +\infty )$ and $q\in (0, 1)\cup (1,
+\infty )$. By analogy with the above right inverse operators
$R_{hs}$ we first construct the operators $\rho_{qs}$, where $s\in
(0,+\infty )$, being the (particular) right inverses of $\delta_q$.
Then, we define the corresponding (particular) right inverses
$P_{qs}$ of $\Delta_q$.
  Namely,
for $q\in (0, 1)$ we have
\begin{equation}\label{sq01rho}
\rho_{qs}f(x)=\left\{
\begin{array}{ccl}
\sum\limits_{m=1}^{-\lfloor log_{q}\frac{s}{x}\rfloor}f(x q^{-m})& & x\in(0, s)\vspace{1mm}\\
0 &\for &x\in [s, sq^{-1})\; ,\\
-\sum\limits_{m=0}^{\lfloor log_{q}\frac{s}{x}\rfloor -1}f(x q^m)& & x\in [sq^{-1}, \infty)\\
\end{array}
\right.
\end{equation}
and for $q\in (1, \infty)$ we have
\begin{equation}\label{sq1infrho}
\rho_{qs}f(x)=\left\{
\begin{array}{ccl}
-\sum\limits_{m=0}^{-\lfloor log_q \frac{x}{s}\rfloor -1}f(x q^{m})& & x\in(0, s)\vspace{1mm}\\
0 &\for &x\in [s, sq)\; .\\
\sum\limits_{m=1}^{\lfloor log_q \frac{x}{s}\rfloor}f(x q^{-m})&& x\in [sq, \infty)\\
\end{array}
\right.
\end{equation}
One can easily verify that $\delta_q\rho_{qs} = I$, for any $s\in
(0,+\infty )$. Now, to find the right inverses $P_{qs}$ of the
divided difference operator $\Delta_q$, defined by (\ref{q-deriv}),
we can write
\begin{equation}\label{decomposeDelta}
\Delta_q = T^{-1}_q\circ\delta_q ,
\end{equation}
where $T_q$ is the invertible operator defined as
\begin{equation}\label{Kq}
T_qf(x)=(q-1)xf(x),
\end{equation}
and apply formula (\ref{composeinverse}), i.e.
$P_{qs}=\rho_{qs}\circ T_q$. For $q\in (0,1)$ the result is
\begin{equation}\label{sq01R}
P_{qs}f(x)=\left\{
\begin{array}{cc}
\sum\limits_{m=1}^{-\lfloor log_q\frac{s}{x}\rfloor }(q-1)xq^{-m}f(x q^{-m})\;&  x\in(0, s)\vspace{1mm}\\
\hspace{2cm}0\hspace{3,3cm} \for &x\in [s, sq^{-1})\; ,\\
-\sum\limits_{m=0}^{\lfloor log_q\frac{s}{x}\rfloor -1}(q-1)xq^mf(x q^m)\; &  x\in [sq^{-1}, \infty)\\
\end{array}
\right.
\end{equation}
and for $q\in (1, \infty)$ there is
\begin{equation}\label{sq1infR}
P_{qs}f(x)=\left\{
\begin{array}{cc}
-\sum\limits_{m=0}^{-\lfloor log_q \frac{x}{s}\rfloor -1}(q-1)xq^{m}f(x q^{m})\;&  x\in(0, s)\vspace{1mm}\\
\hspace{1,8cm}0\hspace{3cm} \for &x\in [s, sq)\; .\\
\sum\limits_{m=1}^{\lfloor log_q \frac{x}{s}\rfloor }(q-1)xq^{-m}f(x q^{-m})\; &  x\in [sq, \infty)\\
\end{array}
\right.
\end{equation}

Although a single right inverse operator can generate all the other
ones by formula (\ref{RA}), the right inverses $P_{qs}$ can be used
to reach certain $q$-antiderivative, the so-called Jackson integral,
being an infinite series, derived formally in \cite{kac}. From this
approach it becomes clear that Jackson integral is not the only
$q$-antiderivative existing and even though it is divergent for
certain function $f$, we can still work with other
$q$-antiderivatives of $f$, well defined by the finite sums, which
are never threatened by a divergency problem.

Namely, in the lower part of formula (\ref{sq01R}) we put
$s\rightarrow 0$ and obtain Jackson integral, being the series
\begin{equation}\label{jackson}
\int f(x)d_qx = (1-q)x\sum\limits_{m=0}^{\infty}q^mf(xq^m),
\end{equation}
for $x\in (0,+\infty )$.

As a next step we formulate
 definite integrals in terms of algebraic analysis and compare
 them with definite $q$-integrals originally defined in $q$-calculus.
 In analogy to formula (\ref{Fa}) let us consider the operators $G_a$
 defined by
 \begin{equation}\label{Ga}
 G_af(x)=f(xq^{-\lfloor log_q\frac{x}{a}\rfloor }),
 \end{equation}
 for any function $f:(0,\infty)\rightarrow\mathbb R$ and $a\in
 (0,\infty )$. Evidently, for $a\in
 (0,\infty )$, operators $G_a$  are surjective onto the family of all q-periodic
 functions defined on $(0,\infty)$. One can also verify the property
 $G_a^2=G_a$, for any $a\in
 (0,\infty )$. Therefore the operators $G_a$ are the initial operators induced by the
 operator $\Delta_q$, for any $a\in
 (0,\infty )$.

Now, according to formula (\ref{definiteI}), we obtain a $q$-definite integral determined  by the initial
operators $G_a$ and $G_b$
 \begin{equation}\label{qG-definiteintegral}
{\cal I}^{b}_{a}=G_b P_{qs} - G_a P_{qs},
 \end{equation}
for any $a,b,s\in (0,+\infty)$ (the above result is independent of
$s$).

In order to interpret formula (\ref{qdefinite}) within this
framework, for any $a, b\in (0,+\infty )$, we should take $q\in (0,
1)$ and sufficiently big positive $s$ for which $a, b\in
[q^{-s+1},+\infty )$, since the last interval corresponds with $(0,
+\infty )$ when $s\rightarrow +\infty$. Assume $a<b=aq^k$, for some
$0>k\in\mathbb Z$ and calculate
$$
{\cal I}^{b}_{a}f(a)=G_b P_{qs}f(a) - G_a P_{qs}f(a)=
P_{qs}f(aq^{-\lfloor log_q\frac{a}{b}\rfloor})-P_{qs}f(aq^{-\lfloor
log_q\frac{a}{a}\rfloor})=
$$
$$=P_{qs}f(aq^{k})-P_{qs}f(a)=
(1-q)aq^{k}\sum\limits_{m=0}^{\lfloor log_q\frac{1}{aq^{k}}\rfloor
-1+s}q^mf(aq^{k}q^m)-$$
$$-(1-q)aq^{k}\sum\limits_{m=0}^{\lfloor log_q\frac{1}{a}\rfloor -1+s}q^mf(a q^m)=
(1-q)a\sum\limits_{m=k}^{-1}q^mf(aq^m).
$$
On the other hand, from formula (\ref{qdefinite}) we obtain
$$
\int\limits_{a}^{b}f(x)d_qx=
\int\limits_{0}^{b}f(x)d_qx-\int\limits_{0}^{a}f(x)d_qx=
$$
$$=(1-q)b\sum\limits_{m=0}^{\infty}q^mf(bq^m)
-(1-q)a\sum\limits_{m=0}^{\infty}q^mf(aq^m)=$$
$$=(1-q)a\sum\limits_{m=k}^{-1}q^mf(aq^m),
$$
which coincides with the previous result. Let us notice that formula
(\ref{qdefinite}) defines $q$-definite integral provided the Jackson
$q$-antiderivative is a convergent series. A simple example of a
function $f$, for which such a formulation of a definite
$q$-integral cannot be applied is $f(x)=\frac{1}{x}\,$, for which
Jackson $q$-antiderivative is evidently divergent. But fortunately,
according to formula (\ref{definiteI}), a definite integral depends
only on the initial operators and is completely independent of a
particular choice of a right inverse used in the calculation.
Therefore, divergency of Jackson integral merely means that this
particular  $q$-antiderivative cannot be used in the calculation of
a given $q$-definite integral.

Let us end this section with the example of a definite $q$-integral
for the above mentioned function $f(x)=\frac{1}{x}\,$, where we
assume $q\in (0, 1)$ and $0<a<b=aq^k$, for some negative
$k\in\mathbb Z$. We obtain
$$
\int\limits_{a}^{b}\frac{1}{x}\,
d_qx=(1-q)a\sum\limits_{m=k}^{-1}q^m\frac{1}{aq^m}=(1-q)\cdot (-k)=
(1-q)log_q\frac{a}{b}\;.
$$
An interesting observation is that the above definite $q$-integral
depends only on the ratio of its limits $a$ and $b$.

\section{Tension spaces}
The usual quantum calculus, i.e. $h$- or $q$-calculus \cite{kac}, is
based on very special difference and divided difference operators.
As one can easily notice, formulae (\ref{h-diff}), (\ref{q-diff})
can be realized for functions defined on an arbitrary set $M$ while
there arises a problem with formulae (\ref{h-deriv}),
(\ref{q-deriv}) since the differences appeared in the corresponding
denominators are undefined unless $M$ is equipped with the usual
algebraic structure. In order to avoid that problem we propose here
to study more general formulation of quantum calculus in a tension
space $(M,\theta )$.

Let $M\neq\emptyset$ and assume the following definition.
\begin{df}\label{theta} By a tension function on $M$ we understand any function
$\theta :M\times M\rightarrow\mathbb R$ such that
\begin{equation}
\theta (p_1,p_2)+\theta (p_2,p_3)=\theta (p_1,p_3) ,
\end{equation}
for any $p_1, p_2, p_3\in M$.
\end{df}
Directly from the above definition, we can prove that any tension
function is skew symmetric, i.e. for any $p_1, p_2\in M$ there is
\begin{equation}
\theta (p_1,p_2) = - \,\theta (p_2,p_1).
\end{equation}
\begin{df}
By a tension space we shall mean a pair $(M,\theta )$, where
$M\neq\emptyset$ and $\theta$ is a tension function on $M$.
\end{df}
In this paper we shall assume that $(M,\theta )$ is a nontrivial
tension space, i.e. there exist points $p, q\in M$ for which
\begin{equation}\label{nontrivial}
\theta (p, q)\neq 0 .
\end{equation}
{\bf Remark:} One can easily check that a linear combination of
tension functions on $M$ is a tension function again. Consequently,
any family $\{\theta^{j}\}_{j\in J}$  of tension functions on $M$
generates the linear space $L=Lin(\{\theta^{j}\}_{j\in J})$, the
so-called tension structure on $M$. Then, by a (multidimensional)
tension space we can understand the pair $(M, L)$. However, in this
paper we consider only a tension space $(M,\theta )$ defined by a
single tension function $\theta$.

With a tension function $\theta$ we shall associate the equivalence
relation in $M$ defined by the formula
\begin{equation}\label{equivelence}
p\sim q\;\;\;\;\ifff\;\;\;\;\theta (p, q)=0.
\end{equation}
Then the equivalence classes of this relation are the following
\begin{equation}\label{classes}
[p\,]=\{q\in M :\; \theta (p,q)=0\}.
\end{equation}

One can easily check that the function $\hat{\theta}$ given by
\begin{equation}\label{quotienttheta}
\hat{\theta}([p\, ], [q\, ])=\theta (p, q) ,
\end{equation}
for $p, q\in M$, is a well defined tension function on the quotient
set $\hat{M}\equiv M/\sim\;$. Thus we have constructed the
"effective" tension space $(\hat{M},\hat{\theta}).$

 On the quotient set $\hat{M}=M/\sim$ we have the natural linear
ordering relation
\begin{equation}\label{linorder}
[p\,]\preceq
[q\,]\;\;\;\;\;\ifff\;\;\;\;\;\hat{\theta}([p\,],[q\,])\leq 0.
\end{equation}
We shall also write $[p\,]\prec [q\,]$ whenever $[p\,]\preceq [q\,]$
and simultaneously $[p\,]\neq [q\,]$.

Then, there is a natural metric $g_\theta$ defined on $\hat{M}$  by
\begin{equation}\label{metric}
g_\theta ([p\,],[q])=|\,\hat{\theta} ([p\,],[q])|\, ,
\end{equation}
for any $p, q\in M$.

In the sequel we will often use mappings
$\theta_q:M\rightarrow\mathbb R$ defined by
\begin{equation}\label{theta_q}
\theta_q (p)=\theta (p,q),
\end{equation}
for any $p, q\in M$. One can easily verify that
$\theta_{q_1}=\theta_{q_2}$, whenever $q_1\sim q_2$.
 Intuitively, the mapping $\theta_q$ we can interpret as a
potential function defined on $M$, associating a scalar potential
$\theta_q(p)$ with any point $p\in M$ and such that $\theta_q(q)=0$
at $q\in M$.

\begin{df}\label{thetadirected}
A mapping $\tau : M\rightarrow M$ is said to be rightward
$\theta$-directed if
\begin{equation}\label{rightward}
[p\,]\prec [\tau (p)\,]\, ,
\end{equation}
 and it is said to be leftward $\theta$-directed if
\begin{equation}\label{leftward}
[\tau (p)\,]\prec [p\,]\, ,
\end{equation}
 for any $p\in M$. We say that $\tau$ is a
 $\theta$-directed mapping if it is either rightward or leftward $\theta$-directed mapping.
\end{df}
Assume the notation: $\tau^0=id_M$ and $\tau^n=\tau\circ\tau^{n-1}$,
for any $n\in\mathbb N\,$.
\begin{p}
For any $\theta$-directed mapping $\tau : M\rightarrow M$ and any
$n\in\mathbb N$, the composition $\tau^n$ has no fixed points, i.e.
\begin{equation}\label{nofixed}
\tau^n(p)\neq p\, ,
\end{equation}
for $p\in M\,$.
\end{p}
{\bf Proof:} Let $\tau $ be a rightward $\theta$-directed mapping.
Then we have inequalities $\theta (\tau (p), p)>0$, ... , $\theta
(\tau^n(p),\tau^{n-1}(p))>0$, for any $n\in\mathbb N$ and $p\in M$.
Consequently,
$$
\theta (\tau^n(p),p)=\theta (\tau^n(p),\tau^{n-1}(p))+\ldots +\theta
(\tau (p), p)>0\, .
$$
Analogously, for a leftward $\theta$-directed mapping we show that
$\theta (\tau^n(p),p)<0\,$, for any $n\in\mathbb N$ and $p\in
M$.\hfill $\Box$

Let us notice that condition (\ref{nofixed}) is not a consequence of
the weaker assumption that $\theta (\tau (p),p)\neq 0$, for any
$p\in M$. In that case there would be $\tau (p)\neq p$ but not
necessarily $\tau^n (p)\neq p\,$, for any $n\in\mathbb N$ and $p\in
M$.
\begin{df}
We say that $\theta$ is homogeneous with respect to $\tau$ (shortly,
$\tau$-homogeneous) if there exists $t\in\mathbb R$, the so-called
$\tau$-homogeneity coefficient, such that
\begin{equation}\label{homogeneous}
\theta (\tau (p_1),\tau (p_2))=t\cdot\theta (p_1,p_2),
\end{equation}
for any $p_1, p_2\in M$.
\end{df}
\begin{p}\label{positivecoeff}Let $\tau :M\rightarrow M$ be a $\theta$-directed mapping
and $\theta$ be a $\tau$-homogeneous tension function. Then, for the
$\tau$-homogeneity coefficient we get
$t>0\, $.
\end{p}
{\bf Proof:} Suppose that $t$ is a $\tau$-homogeneity coefficient
for some $\tau$-homogeneous tension function $\theta$ and assume
that $\tau$ is a $\theta$-directed mapping. Then $\theta
(\tau^2(p),\tau (p))$ and $\theta (\tau (p), p)$ are of common sign
and $\theta (\tau^2(p),\tau (p))=t\cdot\theta (\tau (p), p)$.
Directly from Definition (\ref{thetadirected}) we get $t\neq 0$.
Hence we conclude that $t>0$. \hfill$\Box$
\begin{p}
Let $\theta$ be $\tau$-homogeneous and $\sim$ be the equivalence
relation defined by (\ref{equivelence}). Then, we have the
implication
\begin{equation}\label{implication}
p\sim q\;\;\Rightarrow\;\;\tau (p)\sim\tau (q)\,
\end{equation}
for any $p,q\in M$, or equivalently
\begin{equation}\label{tauinclusion}
\tau ([p\,])\subset [\tau (p)\,],
\end{equation}
for any $p\in M$.
\end{p}
{\bf Proof:} Suppose that $p\sim q$, i.e. $\theta (p, q)=0$. Then we
have $\theta (\tau (p),\tau (q))=\\=t\cdot\theta (p,q)=0$.\hfill
$\Box$

 In general, the inclusion (\ref{tauinclusion})
cannot be inverted, which can be confirmed by the following\\
{\bf Example:} Assume $M={\mathbb R}\times [0, +\infty )$, $\theta
((x_1,y_1), (x_2,y_2))=x_1-x_2$ and $\tau (x,y)=(x+1,y+1)$. Then we
obtain $\tau ([(x,y)])=\{x+1\}\times [1, +\infty )$ and $[\tau
(x,y)]=\{x+1\}\times [0, +\infty )$, i.e. $\tau
([(x,y)])\varsubsetneq [\tau (x,y)]$.
\begin{p}Let $\theta$ be $\tau$-homogeneous with the $\tau$-homogeneity
coefficient $t\neq 0$. Then we have
\begin{equation}
p\nsim q\;\;\;\Rightarrow\;\;\;\tau (p)\nsim\tau (q)\, ,
\end{equation}
for any $p, q\in M$.
\end{p}
 {\bf Proof:}  $\theta (\tau (p), \tau
(q))=t\cdot\theta (p,q)\neq 0$, whenever $p\nsim q$.\hfill$\Box$
\begin{corol}
Let $\theta$ be $\tau$-homogeneous with a $\tau$-homogeneity
coefficient $t$. Assume that $p_0\nsim q_0$ and $\tau (p_0)\sim\tau
(q_0)$, for some $p_0, q_0\in M$. Then $t=0$ and consequently $\tau
(p)\sim\tau (q)$, or equivalently $[\tau (p)\,]=\tau (M)$, for any
$p, q\in M$.
\end{corol}
\section{Quantum ($\tau ,\sigma$)-calculus}\label{tausigmacalculus}
Let  $\sigma , \tau : M \rightarrow M$ be two commuting bijections
and assume ${\cal A}\subset{\mathbb R}^M$ to be a $\sigma^*$,
$\tau^*$-invariant $\mathbb R$-algebra, i.e. $\sigma^*{\cal
A},\tau^*{\cal A}\subset{\cal A}$.
\begin{df}\label{differenceoper}
By the $(\tau ,\sigma )$-quantum differential we mean the mapping
$d_{\tau ,\sigma}: {\cal A}\rightarrow{\cal A}$ given by
\begin{equation}\label{tausigmadifferential}
d_{\tau ,\sigma}f(p)=f(\tau (p))-f(\sigma (p)),
\end{equation}
 for $p\in M$.
\end{df}

One can easily check that the quantum differential $d_{\tau
,\sigma}$ is a linear operator and it fulfills the following Leibniz
product rule
\begin{equation}\label{leibnizform}
d_{\tau ,\sigma}(f\cdot g)(p)=d_{\tau ,\sigma}f(p)\cdot g(\tau
(p))+f(\sigma (p))\cdot d_{\tau ,\sigma}g(p)\, ,
\end{equation}
for any functions $f,g\in{\cal A}$ and $p\in M$.
\begin{df}\label{} By a $(\tau ,\sigma )$-quantum derivation we shall mean any
linear operator $\delta : {\cal A}\rightarrow{\cal A}$ that fulfills
formula (\ref{leibnizform}).
\end{df}
Since the elements $f,g\in{\cal A}$ commute,  the following
combinations are also fulfilled
\begin{equation}\label{ab-combination}
\delta(f\cdot g)(p)= [af(\sigma (p))+bf(\tau(p))]\cdot \delta g(p) +
\delta f(p)\cdot [bg(\sigma (p))+ag(\tau (p))]\, ,
\end{equation}
where $a,b\in\mathbb R$ are coefficients such that $a+b=1$. If
$a\neq b$, formula (\ref{ab-combination}) is equivalent with
(\ref{leibnizform}). In turn, when $a=b=\frac{1}{2}$, formula
(\ref{ab-combination}) becomes symmetric
\begin{equation}\label{ab-sym-combination}
\delta (f\cdot g)(p)= H(f)(p)\cdot \delta g(p) + \delta f(p)\cdot
H(g)(p)\, ,
\end{equation}
where $H(f)(p)=\frac{f(\sigma (p))+f(\tau (p))}{2}\,$. In general,
formula (\ref{ab-sym-combination}) is weaker than
(\ref{leibnizform}) but there exist algebras ${\cal A}$ for which
both formulae are equivalent, i.e. they define the same ${\cal
A}$-module of linear operators (e.g. ${\cal A}={\mathbb R}[x]$,
compare the corresponding comment in Section \ref{hqcalculus}).\\
{\bf Remark:} The mapping $H : {\cal A}\rightarrow {\cal A}$,
defined above, is linear and preserving the unity $1_{\cal A}$ but
in general it is not an algebra homomorphism. The last defect is
precisely the reason why operators defined by
(\ref{ab-sym-combination}) are not differential operators.

Now,  we assume
\begin{equation}\label{theta_neq_0}
[\sigma (p)]_{\tau ,\sigma}\prec [\tau (p)]_{\tau ,\sigma}\, ,
\end{equation}
for any $p\in M$, and define the quantum $(\tau ,\sigma
)$-derivative operator  in a tension space $(M,\theta )$.
\begin{df}\label{quanderivative}
By the  $(\tau ,\sigma )$-quantum derivative we shall mean the
mapping $D_{\tau ,\sigma}: {\cal A}\rightarrow {\cal A}$ given by
\begin{equation}\label{tausigmaderivative}
D_{\tau ,\sigma}f(p)=\frac{d_{\tau ,\sigma}f(p)}{\theta (\tau
(p),\sigma (p))}\equiv\frac{d_{\tau ,\sigma}f(p)}{d_{\tau
,\sigma}\theta_q (p) }\, ,
\end{equation}
for any $f\in{\cal A}$, independently of $q\in M$ .
\end{df}
The assumption (\ref{theta_neq_0}) prevents formula
(\ref{tausigmaderivative}) from zero-valued denominator. However,
owing to the evident symmetry $D_{\tau ,\sigma}=D_{\sigma ,\tau
\,}$, all properties associated with the operator $D_{\tau ,\sigma}$
remain unchanged if the direction of (\ref{theta_neq_0}) is
reversed. Equivalently, relation (\ref{theta_neq_0}) can be
formulated as
\begin{equation}\label{theta_neq_0_1}
[p\,]_{\tau ,\sigma}\prec [\tau \sigma^{-1}(p)]_{\tau ,\sigma}\, ,
\end{equation}
for any $p\in M$. By Definition (\ref{thetadirected}) it means that
$\tau \sigma^{-1}$ is a rightward $\theta$-directed bijection.
Indeed, it is enough to replace $p$ by $\sigma^{-1}(p)$ in formula
(\ref{theta_neq_0}) and obtain (\ref{theta_neq_0_1}).

Evidently,  the quantum derivative
$D_{\tau ,\sigma}$ fulfills the product rule (\ref{leibnizform}).\\

In order to formulate the idea of quantum integration (or the Taylor
interpolation polynomial) we shall need the right inverse operators
defined for the above quantum differential
(\ref{tausigmadifferential}) and quantum derivative
(\ref{tausigmaderivative}).

The following definition will play an important role in our further
analysis.
\begin{df}\label{partition}
We say that a family of subsets $M_k\subset M$, $k\in\mathbb Z$, is
a $(\tau,\sigma )$-partition of $M\neq\emptyset$ if
\begin{itemize}
\item[1)]$\bigcup_{k\in\mathbb Z}M_k=M$,
\item[2)]$M_{k_1}\cap M_{k_2}=\emptyset$, for any $k_1\neq k_2$,
\item[3)]$\tau\sigma^{-1} : M_k\rightarrow M_{k+1}$ is a bijective mapping, for any $k\in\mathbb Z$.
\end{itemize}
\end{df}
 To shorten our notation, the circle symbol "$\circ$"
is omitted for the composition of mappings above and later on.
\begin{p}If $M_k\subset M$, $k\in\mathbb Z$, is a $(\tau,\sigma )$-partition of $M\neq\emptyset$, then
$M_0\neq\emptyset$ and  the composed mapping $(\tau\sigma^{-1})^m$,
for any $m\in\mathbb Z$, has no fixed points.
\end{p}
{\bf Proof:} Suppose $M_0=\emptyset$. Then, by condition (3) we get
$M_k=\emptyset$, for all $k\in\mathbb Z$, which contradicts
condition (1). In turn, let $(\tau\sigma^{-1})^m(p)=p$ for some
$p\in M_k$ and $m\neq 0$. Then by condition (3) we obtain
$p=(\tau\sigma^{-1})^m(p)\in M_k\cap M_{k+m}$ which contradicts
condition (2). \hfill $\Box$\vspace{2mm}

With a given $(\tau,\sigma )$-partition of $M$ we associate the
following integer-valued function $\lfloor\cdot\rfloor_{\tau
,\sigma}: M\rightarrow\mathbb Z$, defined by
\begin{equation}\label{index}
\lfloor p\,\rfloor_{\tau ,\sigma} = k\;\;\; \ifff\;\;\; p\in M_k,
\end{equation}
for any $k\in\mathbb Z$. We shall omit the indices and write
$\lfloor\cdot\rfloor$ whenever $\tau$ and $\sigma$ are fixed.
Automatically, for any $p\in M$, from the above formula we conclude
\begin{equation}
p\in M_{\lfloor p\rfloor}\; .
\end{equation}
\begin{p}For any $p\in M$ there is
\begin{equation}\label{1shift}
\lfloor\tau\sigma^{-1} (p)\rfloor = \lfloor p\, \rfloor +1.
\end{equation}
\end{p}
{\bf Proof.} Let $\lfloor p\,\rfloor = k$, i.e. $p\in M_k$ for some
$k\in\mathbb Z$. Then, $\tau\sigma^{-1}(p)\in M_{k+1}$ and
consequently $\lfloor\tau\sigma^{-1} (p)\rfloor = \lfloor p\,\rfloor
+1$.\hfill $\Box$\vspace{2mm}\\
{\bf Remark:}  Since $\sigma$ is a
bijection, we can always replace $p$ by $\sigma (p)$ and repeat
formula (\ref{1shift}) in the following equivalent form
\begin{equation}\label{1shiftequiv}
\lfloor\tau (p)\rfloor = \lfloor\sigma (p)\rfloor +1.
\end{equation}
\begin{df}\label{partitionfunction}
By a $(\tau ,\sigma )$-partition function (partition function, for
short) of $M$  we mean any integer valued function $\lambda :
M\rightarrow\mathbb Z$ such that
\begin{equation}
\lambda (\tau\sigma^{-1}(p))=\lambda (p)+1\,,
\end{equation}
for any $p\in M$.
\end{df}
One can easily prove the following
\begin{p}
For any $(\tau ,\sigma )$-partition function $\lambda$ of $M$, the
family of sets
\begin{equation}\label{lpartition}
M_k=\lambda^{-1}(k)\subset M\, ,
\end{equation}
where $k\in\mathbb Z$, is a $(\tau ,\sigma )$-partition of $M$.
\end{p}
In the sequel, we say that the $(\tau ,\sigma )$-partition of $M$
given by formula (\ref{lpartition}) is determined by $\lambda$.
Naturally, for a given $(\tau ,\sigma )$-partition of $M$ determined
by $\lambda$ we have
\begin{equation}
\lfloor p\rfloor =\lambda (p)\, ,
\end{equation}
for any $p\in M$. With any $(\tau ,\sigma )$-partition of $M$ we
associate the following
\begin{p}\label{r_tau}
A right inverse of the $(\tau ,\sigma)$-differential $d_{\tau
,\sigma}$ is given by the formula
\begin{equation}\label{rtausigma}
r_{\tau ,\sigma} f(p)=\left\{
\begin{array}{ccl}
-\sum\limits_{m=0}^{-\lfloor p \rfloor -1}f(\tau^m\sigma^{-m-1}(p))& &if \;\;\;\lfloor p \rfloor\leq -1 \\
0&&if\;\;\; \lfloor p\,\rfloor = 0\\
\sum\limits_{m=1}^{\lfloor p\,\rfloor}f(\tau^{-m}\sigma^{m-1}(p))&
&if\;\;\; \lfloor p\,\rfloor\geq 1\; .
\end{array}
\right.
\end{equation}
\end{p}
{\bf Proof:} For $\lfloor\sigma (p)\,\rfloor=k\leq -2$ there is
$\lfloor\tau (p)\,\rfloor=k+1\leq -1$. Then
$$d_{\tau ,\sigma}r_{\tau ,\sigma} f(p)= r_{\tau ,\sigma} f(\tau
(p))-r_{\tau ,\sigma}f(\sigma (p))= -\sum\limits_{m=0}^{-k
-2}f(\tau^{m+1}\sigma^{-(m+1)}(p))+$$
$$+\sum\limits_{m=0}^{-k
-1}f(\tau^{m}\sigma^{-m}(p))= -\sum\limits_{m=1}^{-k
-1}f(\tau^{m}\sigma^{-m}(p))+\sum\limits_{m=0}^{-k
-1}f(\tau^{m}\sigma^{-m}(p))=f(p)\, .$$ For $\lfloor\sigma
(p)\,\rfloor= -1$ there is $\lfloor\tau (p)\,\rfloor=0$. Then
$$d_{\tau
,\sigma}r_{\tau ,\sigma} f(p)= 0-r_{\tau ,\sigma}f(\sigma (p))=
\sum\limits_{m=0}^{-(-1) -1}f(\tau^{m}\sigma^{-m}(p))=f(p).$$ For
$[\sigma (p)\,]_{\tau ,\sigma}=k\geq 1$ there is also $[\tau
(p)\,]_{\tau ,\sigma}=k+1\geq 1$. Then $$d_{\tau ,\sigma}r_{\tau
,\sigma} f(p)= r_{\tau ,\sigma} f(\tau (p))-r_{\tau ,\sigma}f(\sigma
(p))= \sum\limits_{m=1}^{k+1}f(\tau^{-(m-1)}\sigma^{m-1}(p))-
$$
$$-\sum\limits_{m=1}^{k}f(\tau^{-m}\sigma^{m}(p))=
\sum\limits_{m=0}^{k}f(\tau^{-m}\sigma^{m}(p))-
\sum\limits_{m=1}^{k}f(\tau^{-m}\sigma^{m}(p))=f(p).\hspace{5mm}
\Box$$

Next, by using formula (\ref{composeinverse}) we can find the right
inverse $R_{\tau ,\sigma}$ of the $(\tau ,\sigma)$-derivative
$D_{\tau ,\sigma}$.

\begin{p}
A right inverse $R_{\tau ,\sigma}$ of the $(\tau
,\sigma)$-derivative $D_{\tau ,\sigma}$ is given by
\begin{equation}\label{Rtausigma}
 R_{\tau
,\sigma} f(p) =\left\{
\begin{array}{cl}
-\sum\limits_{m=0}^{-\lfloor p\,\rfloor -1}
\theta (\tau^{m+1}\sigma^{-m-1}(p),\tau^m\sigma^{-m}(p))f(\tau^m\sigma^{-m-1}(p))&\iif \;\;\lfloor p\,\rfloor\leq -1 \\
0&\iif\;\; \lfloor p\,\rfloor = 0\\
\sum\limits_{m=1}^{\lfloor p\,\rfloor} \theta
(\tau^{-m+1}\sigma^{m-1}(p),\tau^{-m}\sigma^{m}(p))
f(\tau^{-m}\sigma^{m-1}(p))&\iif\;\; \lfloor p\,\rfloor\geq 1\, .
\end{array}
\right.
\end{equation}
\end{p}
{\bf Proof:} Let us define the operator $T_{\tau ,\sigma}$ by
formula
\begin{equation}\label{Ttausigma}
T_{\tau ,\sigma}f(p)=\theta (\tau (p),\sigma (p))\cdot f(p).
\end{equation}
Thus we write $ D_{\tau ,\sigma}=T_{\tau ,\sigma}^{-1}\circ d_{\tau
,\sigma}$ and using formula (\ref{composeinverse}) we obtain
\begin{equation}\label{Atausigma}
R_{\tau ,\sigma}=r_{\tau ,\sigma}\circ T_{\tau ,\sigma}.
\end{equation}
Finally, we apply (\ref{rtausigma}) and after some calculations
obtain formula (\ref{Rtausigma}).\hfill $\Box$\vspace{2mm}\\
{\bf Remark:} Let us notice that the tension function $\theta$ makes
no explicit contribution on the construction of the right inverse
$r_{\tau ,\sigma}$. The only connection between $r_{\tau ,\sigma}$
and $\theta$ is through formula (\ref{theta_neq_0}) which means that
$\tau\sigma^{-1}$ is a $\theta$-directed mapping. On the other hand,
by formula (\ref{Atausigma}), the right inverse $R_{\tau ,\sigma}$
depends on $\theta$ explicitly.

 Now, let us determine the initial operator $F_{\tau
,\sigma}$ induced by $D_{\tau ,\sigma}\,$ and corresponding with
$R_{\tau ,\sigma}$. Since
\begin{equation}\label{Ftausigma}
F_{\tau ,\sigma}=I-R_{\tau ,\sigma}D_{\tau ,\sigma}=I-r_{\tau
,\sigma}d_{\tau ,\sigma}\, ,
\end{equation}
it becomes  simultaneously the initial operator for $d_{\tau
,\sigma}$ corresponding with $r_{\tau ,\sigma}$.
\begin{p}
The initial operator $F_{\tau ,\sigma}$ induced by $D_{\tau
,\sigma}\,$ and corresponding with $R_{\tau ,\sigma}$  is given by
the formula
\begin{equation}\label{Ftausigma1}
F_{\tau ,\sigma}f(p)=f((\tau\sigma^{-1})^{-\lfloor p\,\rfloor}(p)).
\end{equation}
\end{p}
{\bf Proof:} For $\lfloor p\,\rfloor\leq -1$, we have
$$r_{\tau ,\sigma}d_{\tau ,\sigma}f(p)=
-\sum\limits_{m=0}^{-\lfloor
p\,\rfloor-1}f(\tau\tau^{m}\sigma^{-m-1}(p))+
\sum\limits_{m=0}^{-\lfloor
p\,\rfloor-1}f(\sigma\tau^{m}\sigma^{-m-1}(p))=$$
$$=-\sum\limits_{m=1}^{-\lfloor p\,\rfloor}f(\tau^{m}\sigma^{-m}(p))+\sum\limits_{m=0}^{-\lfloor p\,\rfloor
-1}f(\tau^{m}\sigma^{-m}(p))=f(p)-f((\tau\sigma^{-1})^{-\lfloor
p\,\rfloor}(p))$$ If $\lfloor p\,\rfloor =0$, there is $r_{\tau
,\sigma}d_{\tau ,\sigma}f(p)=0$. For $\lfloor p\,\rfloor\geq 1$, we
have
$$r_{\tau ,\sigma}d_{\tau ,\sigma}f(p)=
\sum\limits_{m=1}^{\lfloor
p\,\rfloor}f(\tau\tau^{-m}\sigma^{m-1}(p))-
\sum\limits_{m=1}^{\lfloor
p\,\rfloor}f(\sigma\tau^{-m}\sigma^{m-1}(p))=$$
$$=\sum\limits_{m=0}^{\lfloor p\,\rfloor} f(\tau^{-m}\sigma^{m}(p))-\sum\limits_{m=1}^{\lfloor p\,\rfloor}
f(\tau^{-m}\sigma^{m}(p))=f(p)-f((\tau\sigma^{-1})^{-\lfloor
p\,\rfloor}(p)). \;\;\;\;\;\;\;\; \Box
 $$

If a $(\tau ,\sigma )$-partition of $M$ is determined by a partition
function $\lambda$, we shall index the right inverses or initial
operators by $\lambda$, i.e. we shall write $r_\lambda\equiv r_{\tau
,\sigma}\,$, $R_\lambda\equiv R_{\tau ,\sigma}$ and $F_\lambda\equiv
F_{\tau ,\sigma}\,$.

 If $\lambda_1$ and $\lambda_2$ are two $(\tau
,\sigma )$-partition functions of $M$ and $R$ is an arbitrary right
inverse of the $(\tau ,\sigma )$-quantum derivative $D_{\tau ,\sigma
}\,$, according to formula (\ref{definiteI}) the corresponding
definite $(\tau ,\sigma )$-integrals are given by
\begin{equation}\label{definiteIII}
{\cal I}^{\lambda_2}_{\lambda_1}=F_{\lambda_2} R - F_{\lambda_1} R\,
.
\end{equation}

{\bf Example:} Let $(M,\theta )$ be a tension space, $D_{\tau
,\sigma}$ be a quantum $(\tau ,\sigma)$-derivation of an algebra
${\cal A}\subset{\mathbb R}^M$ and $\eta$ be another tension
function on $M$ such that the bijective mapping $\tau\sigma^{-1}$ is
$\eta$-directed. Additionally, assume that $\eta$ is $\tau$- and
$\sigma$-homogeneous with both  homogeneity coefficients equal $1$.
Then, for any point $s\in M$, the function $\lambda_{s}$ defined by
\begin{equation}\label{lambdaeta}
\lambda_{s}(p)=\left\lfloor\frac{\eta (p,s)}
 {\eta ({\tau\sigma^{-1}(s),s})}\right\rfloor
\end{equation}
is a $(\tau ,\sigma )$-partition function. In particular, when
$M=\mathbb R$, $\tau (x)=x+h$, $\sigma (x)=x$, $\eta (x,y)=x-y$, for
$x, y, h, s\in\mathbb R$, $h>0$,  we get the partition function
$\lambda_{s}(x)=\lfloor\frac{x-s}{h}\rfloor$ used in $h$-calculus
(see Section \ref{hqcalculus}).
 Hence we obtain the right inverse operators as
well as the  initial operators $F_{\lambda_{s}}$ corresponding with
$\lambda_{s}$. Consequently, the $(\tau ,\sigma )$-definite
integral, for $a, b\in\mathbb R$, is given as
\begin{equation}
{\cal I}_{a}^{b}=F_{\lambda_{b}}R-F_{\lambda_{a}}R\, ,
\end{equation}
where $R$ is an arbitrary right inverse of $D_{\tau ,\sigma}$.

At the end, let us make a comment about higher order $(\tau ,\sigma
)$-difference-like operators. Let $M\neq\emptyset$ and ${\cal
A}_n\subset{\mathbb R}^{M^n}$ be a sequence of $\mathbb R$-algebras,
for $n\in\mathbb N$, and let ${\cal A}={\cal A}_1$. Assume
$(p_1,\ldots ,p_n)\in M^n$ and define $\mu_{p_1,\ldots
,p_n}=\{f\in{\cal A}_n: f(p_1,\ldots ,p_n)=0\}$, the ideal of ${\cal
A}_n$, for any  $n\in\mathbb N$.
\begin{df}
A linear mapping $\Lambda : {\cal A}_1\rightarrow{\cal A}_n$, for a
fixed $n\in\mathbb N$, is said to be of pre-order $n$ if $\Lambda
(\mu_{p_1}\cdot\ldots\cdot\mu_{p_n})\subset\mu_{p_1,\ldots ,p_n}$.
\end{df}
For example, let us explicitely formulate the rule fulfilled by an
operator $\Lambda$ of pre-order $n=1$. From the above definition we
obtain
\begin{equation}
\Lambda ((f_1-f_1(p_1))((f_2-f_2(p_2))\in\mu_{p_1,p_2},
\end{equation}
which means that
\begin{equation}\label{Delta12}
\Lambda ((f_1-f_1(p_1))((f_2-f_2(p_2))(p_1,p_2)=0.
\end{equation}
Formula (\ref{Delta12}) can be written equivalently as
$$
\Lambda (f_1f_2)(p_1,p_2)-f_1(p_1)\Lambda (f_2)(p_1,p_2)
-f_2(p_2)\Lambda (f_1)(p_1,p_2)+
$$
\begin{equation}\label{expand12}
+f_1(p_1)f_2(p_2)\Lambda (1)(p_1,p_2)=0.
\end{equation}
Now, let us define $\delta :{\cal A}\rightarrow{\cal A}$ by formula
\begin{equation}
\delta (f)(p)=\Lambda (f)(\tau (p),\sigma (p))\, ,
\end{equation}
for any $p\in M$. Directly from formula (\ref{expand12}) we obtain
$$
\delta (f_1f_2)(p)- f_2(\tau (p))\delta (f_1)(p) -
  f_1(\sigma (p))\delta (f_2)(p)
+
$$
\begin{equation}\label{expand12delta}
+f_1(\sigma (p))f_2(\tau (p))\delta (1)(p)=0.
\end{equation}
\begin{df}By a quantum $(\tau ,\sigma )$-difference-like operator of order $1$ we shall mean
any operator $\delta_{\tau ,\sigma}$ that fulfills formula
(\ref{expand12delta}). In the case $\delta_{\tau ,\sigma}(1)=0$, an
 operator $\delta_{\tau ,\sigma}$ is
said to be a quantum $(\tau ,\sigma )$-derivative (compare with
formula (\ref{leibnizform})).
\end{df}
In the particular case $\tau =\sigma =id_M$ the above $(\tau ,\sigma
)$-differential operator $\delta_{\tau ,\sigma}$ becomes a usual
differential operator of order $1$ of algebra $\cal A$.

\end{document}